\documentclass[10pt]{amsart}

\usepackage{amsmath}

\usepackage{amssymb}

\usepackage{graphicx}

\newtheorem{thm}{Theorem}[section]
\newtheorem{prop}[thm]{Proposition}
\newtheorem{lemma}[thm]{Lemma}

\theoremstyle{definition}
\newtheorem{defn}[thm]{Definition}

\theoremstyle{remark}
\newtheorem{remark}[thm]{Remark}
\newtheorem{notation}[thm]{Notation}

\numberwithin{equation}{section}

\def\R{\mathbb{R}}

\def\K{\mathbb{K}}
\def\F{\mathcal{F}}

\def\hom{\mathrm{Hom}}
\def\dim{\mathrm{dim\,}}
\def\ker{\mathrm{Ker\,}}
\def\span{\mathrm{span\,}}
\def\A{\mathbf{A}}
\def\g{\mathfrak{g}}
\def\a{\mathfrak{a}}
\def\b{\mathfrak{b}}

\def\IG{\mathfrak{I}}
\def\m{\mathfrak{m}}

\def\p{\mathfrak{p}}

\def\r{\mathfrak{r}}

\begin{document}

\title[A class of overdetermined systems]{
A class of overdetermined systems defined by tableaux:
involutiveness and Cauchy problem}


\author{Emilio Musso}
\address{(E. Musso) Dipartimento di Matematica Pura ed Applicata,
Universit\`a degli Studi di L'Aquila, Via Vetoio, I-67010 Coppito
(L'Aquila), Italy} \email{musso@univaq.it}

\author{Lorenzo Nicolodi}
\address{(L. Nicolodi) Di\-par\-ti\-men\-to di Ma\-te\-ma\-ti\-ca,
Uni\-ver\-si\-t\`a degli Studi di Parma, Parco Area delle Scienze 53/A,
I-43100 Parma, Italy} \email{lorenzo.nicolodi@unipr.it} 

\thanks{Authors partially supported by the MIUR projects
\textit{Propriet\`a geometriche delle variet\`a reali e
complesse} and \textit{Metriche riemanniane e variet\`a differenziali}, and
the GNSAGA of INDAM}


\date{February 28, 2006}

\begin{abstract} 
This article addresses the question of involutiveness and
discusses the initial value problem for a class of overdetermined 
systems of partial differential equations which arise in the theory 
of integrable systems and are defined by tableaux. 
\end{abstract}

\maketitle

\section{Introduction}

The present paper is concerned with  
the overdetermined system of first order nonlinear partial 
differential equations of type
\begin{equation}\label{equation}
Q^a_{\alpha i}\partial_{x^j}F^{\alpha}-Q^a_{\alpha j}\partial_{x^i}F^{\alpha}
=\Phi^a_{ij} \quad  (1\leq i < j \leq n),
\end{equation} 
where $Q^a_{\alpha i}$ ($a = 1,\dots r$; $\alpha = 1, \dots,s$) are constant
coefficients and the non-homogeneous terms
$\Phi_{ij}^a$ are analytic functions depending on both 
the independent and dependent variables. 
The interest in this system originates from the observation that 
many soliton equations of mathematical physics and  
integrable systems occurring in submanifold geometry can be written in the 
form \eqref{equation}. 
Example include 
the equation of harmonic maps from Euclidean or Minkowski $2$-space 
to a Lie group or a symmetric space (the principal chiral model in particle physics) 
\cite{ShaStr,TU, Uhlenbeck};
the Ward equation (or the modified $2+1$ chiral model) \cite{DaiTerng,DaiTerng1};
the curved flat system of Ferus and Pedit \cite{FP1,FP2}; and the
$G/G_0$-system of Terng 
associated to a rank $n$ symmetric space $G/G_0$ \cite{Ter97,BDPT02,TW03},
to name just a few.

\vskip0.2cm

One purpose of this paper is to address the question of involutiveness for the
system \eqref{equation}. 
Roughly speaking, this problem amounts to finding sufficient
conditions under which the system can be solved by repeated use of
equations of Cauchy-Kowalewski type. In the language of exterior differential 
systems (EDS), an involutive system is one that satisfies the 
hypotheses of the Cartan-K\"ahler theorem.
R. Bryant, using the theory of EDS, proved that the $G/G_0$-system 
and the curved flat system are involutive (see Terng and Wang \cite{TW03}). 
This result was taken up and extended in \cite{MNcag} where
the two previous systems were identified as examples of a class of 
involutive systems defined by algebraic tableaux over Lie algebras. 
In this paper, we determine sufficient conditions 
(more general than the involutiveness condition) 
on the 
constant coefficients and the non-homogeneous terms of \eqref{equation}
in order to guarantee that the system becomes involutive after  
a sufficiently large number of successive prolongations.
A second purpose of the work is to discuss the Cauchy problem 
for the system \eqref{equation}, namely formulate
appropriate initial conditions 
and prove a (local) existence and uniqueness theorem for the K-regular analytic 
solutions satisfying given initial conditions. 
As concerns the approach to the problem, these results are in line with those 
obtained by the authors in \cite{MNjmp}. (For another view on this topic
see also \cite{GU}.)

\vskip0.2cm
We use methods by which the above analytic questions may be
reduced to algebraic questions concerning the tableau generated by
the matrices $Q_\alpha = (Q^a_{\alpha i})$ and the associated Spencer complex.
This reduction essentially proceeds in three steps: (1) replacing the consideration of 
the system \eqref{equation} by that of the differential ideal $\IG$ on $M=\R^n\oplus\R^s$
generated by the $2$-forms 
$$
 \Theta^a = Q^a_{\alpha i}dy^{\alpha}\wedge
   dx^i-\Phi^a_{ji}dx^j\wedge dx^i \quad (a=1,\dots,r)
    $$ 
and with independence (transversality) condition 
$\Omega = dx^1\wedge\cdots\wedge dx^n\neq 0$;
(2) replacing the consideration of $(M,\IG)$ by that of its first prolongation 
$(M_{(1)},\IG_{(1)})$,
that is, the restriction to the space of $n$-dimensional integral elements of
$\IG$ of the canonical Pfaffian system on the Grassmann bundle $G_n(TM)$ of 
tangent $n$-planes of $M$;
(3) replacing the consideration of this Pfaffian system and its higher prolongations 
$(M_{(h)},\IG_{(h)})$ by that of the tableau $(Q^a_{\alpha i})$, its prolongations, 
and the harmonic decomposition of the Spencer complex attached to it.


\vskip0.2cm

In Section \ref{s:basic}, we review the needed 
results about the algebraic theory of tableaux, we introduce the Spencer 
complex associated to a tableau, and describe the Gui\-l\-le\-min normal
form for an involutive tableau.

In Section \ref{s:ideals}, we carry out 
the steps of the above list and present
the results concerning the involutiveness of \eqref{equation}.
We prove that if the tableau $(Q^a_{\alpha i})$ is $2$-acyclic and the
non-homogeneous term $\Phi$ satisfies reasonable regularity conditions 
(satisfied  by all the cited examples), then
the prolongation tower
$$
 \cdots \xrightarrow{\pi_{h+1}} (M_{(h)},\IG_{(h)},\Omega)\xrightarrow{\pi_h}
  (M_{(h-1)},\IG_{(h-1)},\Omega) \xrightarrow{\pi_{h-1}}\cdots \xrightarrow{\pi_2}
    (M_{(1)},\IG_{(1)},\Omega)
     $$ 
is indeed a \textit{Frobenius tower}; here we use the terminology introduced by
Bryant and Griffiths in the foundational paper \cite{BG} about the 
characteristic cohomology of differential systems (see Section \ref{s:ideals}
for the definition). 
Next, we prove that if the $k$-th prolongation of the tableau $(Q^a_{\alpha i})$ is 
involutive,
then the $k$-th prolongation $\IG_{(k)}$ of $\IG$ 
is also involutive and has the same characters of the involutive prolongation of 
the tableau. 

In Section \ref{s:examples}, we show that the examples listed above
all share the common formalism \eqref{equation} and fit into the general scheme.

In the final section, we discuss the Cauchy problem 
for the system \eqref{equation}; the existence and uniqueness of K-regular 
analytic solutions
for suitable analytic Cauchy data is proved.
The basic reference on exterior differential systems has been \cite{BCGGG} and
our notation and terminology are consistent with that reference.
Summation convention over repeated indices will be adopted throughout the paper.

\section{Involutive tableaux and the Spencer complex}\label{s:basic}

In this section we briefly review the material and the results we need in the
algebraic theory of tableaux, referring for details to the book by Bryant, et al.
\cite{BCGGG}.


\subsection{Definitions and basic properties}

Let $\a$ and $\b$ be (real or complex) finite dimensional vector spaces. 
A \textit{tableau} $\A$ is a linear subspace of $\hom (\a,\b)$.

An $h$-dimensional subspace $\a_h \subset \a$ is called \textit{generic} with 
respect to $\A$ if the dimension of
$\ker (\A,\a_h) : =\{Q\in \A \,|\, Q_{|_{\a_h}} =0\}$ is a minimum, i.e.,
$$
\dim \ker(\A,\a_h)  =
  \min\{\dim\ker(\A,{\tilde{\a}}_h)\, : \, {\tilde{\a}}_h\in G_h(\a)\},
   $$
where $G_h(\a)$ denotes the Grassmannian of $h$-dimensional subspaces of $\a$.
Similarly, a flag $(0) \subset \a_1\subset \dots \subset \a_n=\a$ of $\a$
is said \textit{generic} if $\a_h$ is generic for all
$h=1,\dots,n$. The set of generic flags is an open and dense subset 
of the flag manifold $\F(\a)$.

The \textit{characters} of $\A$ are the non-negative integers
$s_j(\mathbf{A})$, $j=1,\dots,n$, defined inductively by
$$
   s_1(\A) + \cdots +s_j(\A)= \mathrm{codim\,}
    \ker(\A,\a_j),\quad j = 1,\dots,n,
      $$
where $(0)\subset \a_1\subset\cdots\subset \a_n=\a$ is
a generic flag. 
From the definition it is clear that
$$
 \dim \b \geq s_1\geq
     s_2\geq \dots \geq s_n,\quad \dim\A = s_1+\cdots +s_n.
     $$ 
If $s_{\nu}\neq 0$, but $s_{\nu +1}=0$, we say that $\A$ has 
\textit{principal character} $s_{\nu}$ and call $\nu$ the 
\textit{Cartan integer} of $\A$.


The \textit{first prolongation} of $\A$ is the subspace
${\A}^{(1)}\subset \A \otimes \a^\ast \simeq \hom(\a,\A)$ 
of all linear transformations $Q :\a\to \A$ such that
$$
Q(X)(Y) = Q(Y)(X), \quad\text{for all\,}\, X,Y\in \a.
$$
Notice that ${\A}^{(1)}$ is itself a tableau and that 
$
{\A}^{(1)} = \left(\A\otimes \a^\ast\right) \cap \left(\b\otimes S^2(\a^\ast)\right)
  $.
(By $S^h(\a^\ast)$ we denote the symmetric $h$-fold tensor product of $\a^\ast$.)
The \textit{$h$-th prolongation} of ${\A}$ is defined
inductively by setting ${\A}^{(h)} = {{\A}^{(h-1)}}^{(1)}$, 
for $h\geq 1$, where ${\A}^{(0)}= \A$ and by convention 
${\A}^{(-1)}= \b$. 
The prolongation ${\A}^{(h)}$ is identified with
$$
  \mathbf{A}^{(h)}=\left(\mathbf{A}\otimes S^h(\a^*)\right)\cap 
     \left(\b\otimes S^{h+1}(\a^\ast)\right)
       $$
and an element $Q_{(h)}\in \b\otimes S^{h+1}(\a^\ast)$ belongs to ${\A}^{(h)}$ 
if and only if $i(X)Q_{(h)} \in {\A}^{(h-1)}$, for all $X \in \a$.
The direct sum
$$
 \A^{[\infty]}= \bigoplus_{h\ge 0}{\A}^{(h)},
  $$
is called \textit{total prolongation} of $\A$.
The total prolongation has a natural filtration
$$
 \A={\A}^{[0]}\subset {\A}^{[1]}\subset {\A}^{[2]}\subset \cdots
  \subset
   {\A}^{[p]}\subset \cdots\subset {\A}^{(\infty)},
    $$ 
where ${\A}^{[p]}= \bigoplus_{0\leq h\leq p}{\A}^{(h)}$.


\vskip0.3cm

One can establish the following result (cf. \cite{BCGGG}) 
\begin{equation}\label{inequality}
 \dim {\A}^{(1)} \leq s_1 +2s_2 + \cdots + ns_n.
  \end{equation}
$\A$ is said \textit{involutive} (or \textit{in involution}) if equality holds 
in the above inequality.

Also, \textit{for any tableau $\A$, there exists an integer $h_0$ such that 
${\A}^{(h)}$ is involutive, for all $h\geq h_0$}. 
We call the smallest integer $k$ with this property
the \textit{involutive index} and the 
corresponding prolongation the \textit{involutive 
prolongation} of the tableau. The \textit{involutive characters} of $\A$ are
the characters of the involutive prolongation ${\A}^{(k)}$. 
Accordingly, we have the notion of \textit{involutive principal character}
and \textit{involutive Cartan integer} of a tableau.

Another basic property is that \textit{every prolongation of an involutive 
tableau is involutive}. In particular,
if $\A$ is involutive, 
the characters of $\A$ and ${\A}^{(1)}$ are related by
$$
  s^{(1)}_j: = s_j({\A}^{(1)}) = s_n(\A)+\cdots + s_j(\A) \qquad (j=1,\dots,n).
     $$ 
This implies that 
\textit{the principal character and
the Cartan integer are invariant under prolongation of an involutive tableau}.

\subsection{The Spencer complex}
With the previous notation, let
$$
  \delta^{q,p} : \b\otimes S^{q}(\a^\ast) \otimes \Lambda^p(\a^\ast) \rightarrow  
    \b\otimes S^{q-1}(\a^\ast) \otimes \Lambda^{p+1}(\a^\ast)
  $$
be the operator given by
\begin{equation}\label{delta}
 \delta^{q,p}\phi(X_1,\dots,X_{p+1}) : = 
  \sum_k (-1)^{k+1}\,i(X_k)\phi(X_1,\dots,\hat{X}_k,\dots,X_{p+1}).
   \end{equation}
If $q=0$ we set $\delta^{0,p} =0$, for $p\geq 0$. From the definition, 
we have that $\delta\circ\delta = 0$. Also, 
the sequence of the corresponding bi-graded complex is exact except when $q=0$
and $p=0$. 

\vskip0.2cm
Let $\A\subset \hom(\a,\b)$ be a tableau with
prolongations ${\A}^{(h)}$, $h\geq 0$.
Consider the sequence of spaces 
$$
 C^{q,p}(\A):={\A}^{(q-1)}\otimes \Lambda^p(\a^\ast),
  $$
for integers $q\geq 0$ and $0 \le p \leq n$.
Since ${\A}^{(q-1)} \subset \b \otimes S^{q}(\a^\ast)$, the space
$C^{q,p}(\mathbf{A})$ is a subspace of 
$\b\otimes S^{q}(\a^\ast) \otimes \Lambda^p(\a^\ast)$.
Now the sequence of subspaces $C^{q,p}(\A)$ is stable under 
$\delta$, that is,
$\delta C^{q,p}(\A) \subset C^{q-1,p+1}(\A)$.
In fact, for any $\phi \in C^{q,p}(\A)$, each summand in the r.h.s.
of \eqref{delta} 
lies in ${\A}^{(q-2)}$ by the defining property of prolongations.
We still have $\delta\circ\delta =0$ (for $q\geq 1$ and $p=1,\dots,n$), but 
the sequence 
$$
 C^{q+1,p-1}(\A) \xrightarrow{\delta^{q+1,p-1}} C^{q,p}(\A) 
\xrightarrow{\delta^{q,p}} C^{q-1,p+1}(\A) 
  $$
is no longer exact for all $p$ and $q$.
The associated cohomology groups 
$$H^{q,p}(\A) := Z^{q,p}(\A)/B^{q,p}(\A)$$
are called the \textit{Spencer groups} of $\A$,
where $B^{q,p}(\A)=\mathrm{Im}(\delta^{q+1,p-1})$ and
$Z^{q,p}(\A)=\ker(\delta^{q,p})$.
Notice that
$Z^{0,p}(\A)=\b\otimes\Lambda^p(\a^\ast)$ and 
   $Z^{q,1}(\A)={\A}^{(q)}$, for all $q\geq 1$, $p\geq 0$.

\begin{remark}
A decomposable element of $C^{q,p}(\A)$, $q\ge 1$, is of the form 
$
  \mathcal{Q}=Q\otimes
     (\alpha^{i_1}\wedge\cdots\wedge \alpha^{i_p})
    $,
where $Q\in{\A}^{(q-1)}\subset \hom(\a,{\A}^{(q-2)})$. 
If $(\overline{Q}^{(q-2)}_1,\dots,\overline{Q}^{(q-2)}_{r_{q-2}})$ is a 
basis of ${\A}^{(q-2)}$, then $Q=Q^a_i \overline{Q}^{(q-2)}_{a}\otimes \eta^i$ 
and
$$
  Q\otimes (\alpha^{i_1}\wedge\cdots\wedge \alpha^{i_p}) =
    Q^a_i(\overline{Q}^{(q-2)}_{a}\otimes \eta^i)\otimes
     (\alpha^{i_1}\wedge\cdots\wedge \alpha^{i_p}).
      $$ 
Then $\delta(\mathcal{Q})$ is the element of $C^{q-1,p+1}(\A)$ defined by
$$
 \delta(\mathcal{Q})=Q^a_i
  \overline{Q}^{(q-2)}_{a}\otimes (\eta^i\wedge
   \alpha^{i_1}\wedge\cdots\wedge \alpha^{i_p}).
      $$ 
\end{remark}


\vskip0.2cm
A significant result in the subject is that the vanishing of the $H^{q,p}$ is equivalent 
to involutiveness.

\begin{thm}
A tableau $\A$ is involutive if and only if $H^{q,p}(\A)$ is zero, 
for all $q\geq 1$ and $p\geq 0$.
\end{thm}

If $H^{q,2}(\A)=(0)$, for all $q\ge 1$, the tableau $\A$ is called 
\textit{$2$-acyclic}. 

\vskip0.3cm

Suppose that $\a$ and $\b$ are endowed with (positive definite) scalar 
products, and consider the scalar product induced on 
$C^{q,p}(\A)={\A}^{(p-1)} \otimes\Lambda^p(\a^\ast)$.
With respect to this scalar product, let
$$
 \delta^\ast : C^{q-1,p+1}(\A)\to C^{q,p}(\A)
     $$ 
be the \textit{co-differential} defined by 
$\langle \delta^\ast(\zeta),\rho\rangle = \langle \zeta, \delta(\rho)\rangle$, 
for $\zeta\in C^{q-1,p+1}(\A)$ and $\rho\in C^{q,p}(\A)$. 
This yields a \textit{harmonic decomposition} of the Spencer complex:
$$
 C^{q,p}(\A) = B^{q,p}(\A)\oplus
   \mathcal{H}^{q,p}(\A)\oplus B_{q,p}(\A),
   $$
where $B^{q,p}(\A) = \delta (C^{q+1,p-1}(\A))$,
$B_{q,p}(\A) =\delta^\ast(C^{q-1,p+1}(\A))$, and
$\mathcal{H}^{q,p}(\A) =\ker(\delta)\cap \ker(\delta^\ast)$.
   From the definition, it follows that:
\begin{itemize}
\item the application
$\zeta \in \mathcal{H}^{q,p}(\A)\to [\zeta]\in H^{q,p}(\A)$ is an isomorphism;

\item 
$\ker(\delta) = Z^{q,p}(\A)= B^{q,p}(\A)\oplus \mathcal{H}^{q,p}(\A)$;

\item 
$\ker(\delta^\ast) = Z_{q,p}(\A)= \mathcal{H}^{q,p}(\A)\oplus B_{q,p}(\A)$;
     
\item 
$\mathcal{H}^{q,p}(\A)=Z^{q,p}(\A)\cap Z_{q,p}(\A)$. 

\end{itemize}

\noindent Finally, observe that the restriction of the Spencer differential 
to $B_{q,p}(\A)$ is an isomorphism onto
$B^{q-1,p+1}(\A)$. The inverse of
$\delta_{|_{B_{q,p}}(\A)}$ will be denoted by
$$
  \sigma_{q,p}: B^{q-1,p+1}(\A)\to B_{q,p}(\A).
    $$

\subsection{The Guillemin normal form for involutive tableaux}\label{ss:guillemin}

Let $\A\subset \hom(\a,\b)$ be a tableau.
Let $\mathcal{A}=(A_1,\dots,A_n)$ denote any basis of
$\a$, $(\alpha^1,\dots,\alpha^n)$ its dual basis, and
$\mathcal{B}=(B_1,\dots,B_r)$ any basis of $\b$.
Identify $\hom(\a,\b)$ with the space of $r\times n$ matrices $\K(r,n)$
and consider $\A$ as a subspace of $\K(r,n)$. Let $\pi = (\pi^a_j)$ be 
the tautological $1$-form 
on $\hom(\a,\b)$ and let $\pi_1,\dots,\pi_n$ denote its column vectors. 
Set $\pi_j$ $=$ $(\pi_j^{[\nu]},\dots,\pi_j^{[0]})^T$, where 
$\pi_j^{[\nu]}$ $=$ $(\pi_j^{1},\dots,\pi_j^{s_{\nu}})^T$,
$\pi_j^{[\rho]}$ $=$ $(\pi_j^{s_{\nu}+\cdots+s_{\rho+1}+1}$, $\dots$,
$\pi_j^{s_{\nu}+\cdots+s_{\rho}})^T$, for
$\rho=\nu-1,\dots,1$, and
$\pi_j^{[0]}$ $=$ $(\pi_j^{s_{\nu}+\cdots + s_{1}+1},\dots,\pi_j^r)^T$.


\vskip0.2cm
The following result is essentially due to Guillemin (cf. \cite{BCGGG}, \cite{Gui68}).

\begin{prop}
If $\A \subset \hom(\a,\b)$ is involutive, the bases $\mathcal{A}$ and 
$\mathcal{B}$ can be chosen so that 
$\pi_i^{[0]}|_{\A} = 0$, for $\,i\ge 1$,
$\pi_i^{[\rho]}|_{\A}\in\span(\pi_1^1|_{\A},\dots,\pi_1^{s_1}|_{\A}$, $\dots,
\pi_{\rho}^1|_{\A},\dots$, $\pi_{\rho}^{s_{\rho}}|_{\A})$, for $\rho=1,\dots,\nu$ 
and $i \geq \rho$. In particular,
$(\pi_1^1|_{\A},\dots,\pi_1^{s_1}|_{\A}$, $\pi_2^1|_{\A},\dots,\pi_2^{s_2}|_{\A}$,
$\dots,\pi_{\nu}^1|_{\A}$, $\dots,\pi_{\nu}^{s_{\nu}}|_{\A})$ is a basis of
${\A}^*$ and the elements of its dual basis
$$ 
 \mathcal{Q}=(Q_{[1],1},\dots,Q_{[1],s_1},Q_{[2],1},\dots,Q_{[2],s_2},\dots,
   Q_{[\nu],1},\dots,Q_{[\nu],s_{\nu}})
    $$
can be expressed by 
\begin{equation}\label{guillemin'snormalform}
\begin{array}{c}
 Q_{[j],a} = B_a\otimes \alpha^j+\sum_{h=j+1}^{\nu}\left(
   \sum_{b=s_h+1}^{s_{h+1}}Q_{[j],a,h}^{\,\,b}B_b\otimes
    \alpha^h\right) +\\
     +\sum_{h=\nu+1}^{n}\left(\sum_{b=1}^{s_{\nu}}
       Q_{[j],a,h}^{\,\,b}B_b\otimes \alpha^h\right),
     \end{array}
        \end{equation}
for all $j=1,\dots,\nu$; $a=1,\dots,s_j$.
\end{prop}

We say that $\mathcal{Q}$ is a \textit{normal basis} of the involutive tableau $\A$
and that the expression \eqref{guillemin'snormalform} is the 
\textit{Guillemin normal form} for the elements of $\mathcal{Q}$ with respect to the
basis $\mathcal{A}$ and $\mathcal{B}$.

\begin{remark}
Notice that the basis $(A_1,\dots,A_n)$ of $\a$ is automatically
\textit{generic} with respect to the tableau $\A$, that is, the corresponding flag:
$(0)\subset \a_1=\span(A_1)\subset
    \a_2=\span(A_1,A_2)\subset \cdots\subset
       \a=\span(A_1,\dots,A_n)$  
is generic.
\end{remark}

\section{Differential ideals associated with tableaux}\label{s:ideals}
 
\begin{notation}
We start by recalling some terminology and notation.
An exterior differential system (EDS) is a pair $(M,\mathfrak{I})$ consisting 
of a smooth manifold $M$ and a differential ideal $\mathfrak{I}$, i.e., a homogeneous, 
differentially closed ideal 
$\IG \subseteq \Omega^\ast(M)$ in the algebra of smooth differential forms 
on $M$. A Pfaffian differential system is an EDS which is 
differentially generated by $1$-forms. It is given by the sections of a smooth 
subbundle $I\subset T^\ast M$ of the cotangent bundle; 
the space of section of this subbundle, and sometimes the subbundle itself, 
is also referred to as a Pfaffian system. We use the
notation $\{\alpha, \beta, \dots\}$ for the (two-sided) algebraic ideal
generated by forms $\alpha, \beta, \dots$, and $\{I\}$ for the algebraic ideal
generated by the sections of a Pfaffian system $I\subseteq T^\ast M$.
\end{notation}

With reference to the system \eqref{equation}, let $\a =\R^n$, $\b=\R^{s}$ 
and denote by $\A\subset \hom(\a,\b)$ the tableau
generated by the matrices $Q_{\alpha}$ $=$ $(Q_{\alpha i}^a)\in
\K(r,n)$, $\alpha = 1,\dots,s$.
Assume that the $Q_{\alpha}$
be linearly independent so that $\b$ can be identified with $\A$, 
and let $q^{\alpha}_{(0)}$ denote the coordinates 
with respect to $Q_{\alpha}$.
Let $\IG \subset \Omega^*(M)$ be the differential ideal on $M=\a\oplus\A$ 
generated by the $2$-forms 
$$
  \Theta^a=Q^a_{\alpha j}dq^{\alpha}_{(0)}\wedge dx^j
    -\Phi_{ij}^adx^i\wedge dx^j \quad (a=1,\dots,r),
      $$
where $\Phi_{ij}^a$ are analytic maps of the $x^i$ and the $q^{\alpha}_{(0)}$, and 
$\Omega=dx^1\wedge\cdots\wedge dx^n\neq 0$ gives the independence (transversality) 
condition.
Locally, the integral manifolds of $(\IG,\Omega)$ can be expressed as graphs
of solutions 
to \eqref{equation}. Further, we assume that 
$\Phi(x,q) = (\Phi_{ij}^a) : M \to C^{0,2}(\A)$ 
takes values in $B^{0,2}(\A)$ and 
satisfies
\begin{equation}\label{torsionezero}
 \sum_{(1,2,3)}\Phi_\ast|_{(A,Q)}(A_1+Q_{(1)}(A_1))(A_2,A_3)=0
  \end{equation} 
(sum over cyclic permutations),
for $(A,Q)\in M$, $A_1,A_2,A_3\in \a$, and $Q_{(1)}\in {\A}^{(1)}$, where 
$\Phi_\ast : T_{(A,Q)}M\simeq M \to T_{\Phi(A,Q)}B^{0,2}(\A)\simeq B^{0,2}(\A)$ 
denotes the differential of $\Phi$ at $(A,Q)$.
(In the next section, we will see that the examples cited in the introduction
all satisfy these conditions.)
 
\begin{defn}
We call $(M, \IG,\Omega)$ the \textit{differential system associated with
the tab\-leau $\A$ and the non-homogeneous term $\Phi$}. 
\end{defn}
 
Next, consider a harmonic decomposition of the Spencer complex, and in particular
the decomposition
$$
 C^{1,1}(\A)=B^{1,1}(\A)\oplus
   \mathcal{H}^{1,1}(\A)\oplus B_{1,1}(\A).
      $$
Let $S_{(1)} : M\to B_{1,1}(\A)$ be the analytic mapping defined by 
$\delta^{1,1}(S_{(1)})=\Phi$.
Then 
$$
  \Theta^a = Q^a_{\alpha j}(dq^{\alpha}_{(0)}+{S_{(1)}}_i^{\alpha}dx^i)\wedge dx^j,
    $$ 
from which follows that the $n$-dimensional integral elements
of the differential system $(\IG,\Omega)$ are generated by vectors
of the form 
$$
 v_i=\partial_{x^i}-({S_{(1)}}_i^{\alpha}+{Q_{(1)}}_i^{\alpha})\partial_{q^{\alpha}_{(0)}}
  \quad (i=1,\dots,n),
     $$
where $Q_{(1)}=({Q_{(1)}}_j^{\alpha})\in {\A}^{(1)}$. 
Therefore, the space $V_n(\IG,\Omega)$ of $n$-dimensional integral elements 
can be identified with $M_{(1)}=\a\oplus {\A}^{[1]}$ and the restriction 
to $V_n(\IG,\Omega)$ of the canonical
contact Pfaffian system on the Grassmann bundle $G_n(TM)$ 
(i.e., the first prolongation of $\IG$) 
coincides with the Pfaffian system given by the subbundle
$I_{(1)}\subset T^\ast(M_{(1)})$,
whose sections are the components of the $1$-form
$$
 \beta_{(0)}: = dQ_{(0)}+(S_{(1)}+Q_{(1)})\cdot dx,
   $$ so that locally
$I_{(1)} = \span\{\beta_{(0)}\}$.

\begin{defn}
Let $(M_{(k)}, \IG_{(k)}, \Omega)$, $k\geq 1$, be a sequence of
smooth manifolds $M_{(k)}$ and Pfaffian systems $\IG_{(k)}$
with surjective submersions $\pi_k : M_{(k)} \to M_{(k-1)}$ such that 
$\pi_k^\ast (\IG_{(k-1)})
\subseteq \IG_{(k)}$. Let $I_{(k)}\subset T^\ast M_{(k)}$ be the
subbundle whose sections are the $1$-forms in $\IG_{(k)}$. The above data
define a \textit{Frobenius tower} if $d I_{(k)} \equiv 0 \mod \{I_{(k+1)}\}$
(cf. \cite{BG}).
\end{defn}
 
We are now in the position to prove the following. 

\begin{thm}\label{theoremA}
Let $\A \subset \hom(\a,\b)$ be $2$-acyclic and assume $\Phi : M \to B^{0,2}(\A)$ 
is analytic and satisfies \eqref{torsionezero}. Then the following statements hold true:

\begin{enumerate}

\item $(M_{(1)},\IG_{(1)},\Omega)$ is a quasi-linear Pfaffian system 
       with vanishing torsion;

\item  $(M_{(1)},\IG_{(1)},\Omega)$ admits a \textsl{Frobenius tower} 
$$
\cdots \rightarrow (M_{(h)},\IG_{(h)},\Omega)\rightarrow
  (M_{(h-1)},\IG_{(h-1)},\Omega)\rightarrow\cdots \rightarrow
    (M_{(1)},\IG_{(1)},\Omega);
      $$

\item 
if $k$ is the least integer such that ${\A}^{(k)}$ is involutive,
then $(M_{(k)},\IG_{(k)},\Omega)$ is involutive and its Cartan characters 
coincide with the characters of ${\A}^{(k)}$.
\end{enumerate}
\end{thm}
 
\begin{proof}
Let\footnote{If $\phi$ and $\psi$
are matrices of real-valued differential forms, $\phi \dot{\wedge} \psi$ 
denotes the matrix multiplication where the entries are multiplied via 
wedge product.} 
$$
 \overline{\mathcal{D}}S_{(1)}:=(\partial_x S_{(1)}+\partial_{Q_{(0)}}
   S_{(1)}(S_{(1)}+Q_{(1)})(dx))\dot{\wedge}dx.
    $$  
Then $\Phi$ satisfies \eqref{torsionezero} if and only if 
$\delta^{1,2}(\overline{\mathcal{D}}S_{(1)})=0$. 
Since $\mathbf{A}$
is $2$-acyclic, there exists a unique analytic map 
$S_{(2)}: M_{(1)}\to B_{2,1}(\A)$ such that
$\delta^{2,1}(S_{(2)}) = -\overline{\mathcal{D}}S_{(1)}$. 
With this notation, the structure equations of $(\IG_{(1)},\Omega)$
become
$$
 d\beta_{(0)}\equiv -(dQ_{(1)}-S_{(2)}(dx))
  \dot{\wedge}dx \mod\{I_{(1)}\}.
   $$ 
This yields that
$(\IG_{(1)},\Omega)$ is a linear Pfaffian system with vanishing torsion,
proving $(1)$.

\vskip0,2cm
For result $(2)$, we proceed by induction on $h$. For each integer $h\ge 1$, let
$M_{(h)} := {\a}\oplus {\A}^{[h]}$ and denote its elements
by $(x,Q_{(0)},\dots,Q_{(h)})$.
Consider the analytic map $S_{(h+1)} : M_{(h)}\to
B_{h+1,1}(\A)$ defined inductively by
\[ 
\delta^{1,1}(S_{(1)})  = \Phi;\qquad
\delta^{r,1}(S_{(r)})  = -\overline{\mathcal{D}}(S_{(r-1)}) \quad (r=2,\dots,h),
\]
where
\begin{equation}\label{totalderivation}
 \overline{\mathcal{D}}(S_{(r)})=((\partial_x
   S_{(r)}+\sum_{s=0}^{r-1}\partial_{Q_{(s)}}
     S_{(r)}(S_{(s+1)}+Q_{(s+1)}))\dot{\wedge}dx )\dot{\wedge}dx.
      \end{equation}
Consider the $1$-forms
$$
 \beta_{(r)}=dQ_{(r)}-(S_{(r+1)}+Q_{(r+1)})(dx) \quad (r=0,\dots,h-1),
  $$
$$
  \pi_{(h)}=dQ_{(h)}-S_{(h+1)}(dx),
    $$ 
and the filtration of subbundles locally given by
$$
 I_{(1)}=\span\{\beta_{(0)}\}\subset\cdots\subset
   I_{(h)}=\span\{\beta_{(0)},\cdots,\beta_{(h-1)}\}\subset
     T^*M_{(h)}
      $$ 
with independence condition $\Omega=dx^1\wedge\cdots\wedge dx^n$.
A direct computation yields that the \textit{structure equations} of the Pfaffian 
system $(\IG_{(h)},\Omega)$ can be written in the form:
\[
 \begin{aligned}
  d\beta_{(0)} &\equiv -\beta_{(1)}\dot{\wedge}dx,\mod \{I_{(1)}\},\\
   d\beta_{(1)} &\equiv -\beta_{(2)}\dot{\wedge}dx,\mod \{I_{(2)}\},\\
    & \cdots \cdots \\
      d\beta_{(h-2)} &\equiv -\beta_{(h-1)}\dot{\wedge}dx \mod \{I_{(h-1)}\},\\
       d\beta_{(h-1)} &\equiv -\pi_{(h)}\dot{\wedge}dx \mod \{I_{(h)}\}.
        \end{aligned}
         \]
This implies that the prolongation tower
\[
 \cdots \rightarrow(M_{(h)},\IG_{(h)},\Omega)
  \rightarrow (M_{(h-1)},\IG_{(h-1)},\Omega)
   \rightarrow \cdots\rightarrow (M_{(1)},\IG_{(1)},\Omega)
    \] 
is indeed a Frobenius tower (cf. \cite{BG}). The general inductive step is clear.

\vskip0.2cm

As for $(3)$, let $k$ be the least integer such that $\A^{(k)}$ is involutive and
let $s_1^{(k)}$, $\dots$, $s_n^{(k)}$ be the involutive characters. Consider the
Pfaffian system $(M_{(k)},\IG_{(k)},\Omega)$. From the structure equations,
we see that $\pi_{(h)}$ is the \textit{tableau matrix} of the system.
This means that the \textit{reduced tableau matrix} is the tautological 
form of ${\A}^{(k)}$. Knowing that 
${\A}^{(k)}$ is involutive and using
the Cartan test for involution of linear Pfaffian systems, it follows that
$(\IG_{(k)},\Omega)$ is involutive, with Cartan characters $s_1^{(k)}$, $\dots$, 
$s_n^{(k)}$.

\end{proof}

\begin{remark} 
Let $M_{(\infty)} :=\a\oplus {\A}^{[\infty]}$ and let
$F:M_{(\infty)}\to \K$ be a function such that $F=\pi_{h}\circ
f$, where $\pi_h : M_{(\infty)}\to M_{(h)}$ is the projection and
$f:M_{(h)}\to \K$ is a differentiable (resp.,
analytic, holomorphic) function. Then $F$ is said of class $C^{\infty}$
(resp., analytic, holomorphic). Let $\Omega^*(M_{(\infty)})$ be the
graded differential algebra
${C}^{\infty}(M_{(\infty)},\K)\otimes_{\K}(\oplus_{p=0}^{\infty}
\Lambda^p (M_{(\infty)}))$ 
(the exterior differentiation is defined as in the finite-dimensional
case) and let
$\IG^{\infty}\subset \Omega^*(M_{(\infty)})$ be the Frobenius ideal
generated by $\{\beta_{(h)}\}_{h\in \mathbb{N}}$. 
The \textit{characterisic complex}\footnote{For a detailed description
of the characteristic complex and of the corresponding cohomology theory we
refer the reader to the article of Bryant and Griffiths \cite{BG}.} of
the differential system $(M,\IG,\Omega)$ is the quotient
$\overline{\Omega}=\Omega^*(M_{(\infty)})/\IG^{\infty}$
endowed with the coboundary operator induced, passing to quotient,
by the exterior differentiation of $\Omega^*(M_{(\infty)})$. 
In our case
$\overline{\Omega} = {C}^{\infty}(M_{(\infty)},\K)\otimes_{\K}\Lambda^*(\a)$
and 
$$
 \overline{d}(F_{i_1\dots{i_q}}dx^{i_1}\wedge\cdots\wedge
   dx^{i_q})=\left(\mathcal{D}(F_{i_1\dots i_q})\cdot dx\right)\wedge
    dx^{i_1}\wedge\cdots\wedge dx^{i_q},
     $$ 
where
$\mathcal{D}(F)=\partial_x F +\sum_{h\ge 0}\partial_{Q_{(h)}}
  F(S_{(h+1)}+Q_{(h+1)})$ 
is the operator of \textit{total derivation}.
The computation of the $S_{(r)}$, though complicated, is essentially 
a problem of linear algebra. This suggests the possibility of
explicitly computing
the characteristic cohomology of the Pfaffian systems considered 
in our discussion; 
details of this computation will appear in a subsequent publication.
\end{remark}

\section{Examples}\label{s:examples}

\subsection{Systems associated with a Cartan decomposition} 

Let $\g=\g_0\oplus \m$ be a 
Cartan decomposition of a real semisimple Lie algebra $\g$. We have the relations
$[\g_0,\g_0]\subset \g_0$, $[\g_0,\m]\subset \m$, $[\m,\m]\subset
\g_0$, and $\g_0$ and $\m$ are orthogonal with respect to the Killing form $B$
of $\g$. 
Moreover, the restrictions of $B$ to $\g_0$ and $\m$
are negative and positive definite, respectively. 

Let $\a\subset \m$ be a maximal ($n$-dimensional) abelian subspace. 
The subspace $\m$ decomposes as $\m=\a\oplus \b$, where
$\b = \a^\perp$ is the orthogonal complement with respect to $B$. Next,
let
$$ 
  \g_{\a}=\{X\in \g_0 \,:\,[X,\a]=0\}.
     $$ 
Then $\g_0= \g_{\a}\oplus \p$, where $\p$ denotes the orthogonal complement of 
$\g_{\a}$ in $\g_0$ with respect to the Killing form.
Now, if $A\in \a$ is a regular element, the mappings 
$$
 \text{ad}_A|_{\b}:\b\to \p, \quad \text{ad}_A|_{\p}:\p\to \b
   $$ 
define vector space isomorphisms between $\b$ and $\p$. 
Moreover, since $[\p,\a]\subset \b$ and the map 
$$
 \p \ni X\mapsto -\text{ad}_X|_{\a}\in \hom(\a,\b)
     $$ 
is injective, $\p$ can be realized as a tableau, say $\mathbf{P}$, in $\hom(\a,\b)$.
Similarly, as $[\b, \a]\subset \p$ and the map 
$\b \ni B \mapsto \text{ad}_B|_{\a}\in \hom(\a,\p)$
is injective, $\b$ can be realized as a tableau $\A$ in $\hom(\a,\p)$.

Now let $\Phi : \A\to \Lambda^2(\a^\ast)\otimes \p$ be the polynomial map
defined by 
\[
  \Phi|_B(A_1,A_2) : =[[A_1,B],[A_2,B]].
   \]
$\A$ and $\Phi$ define the following system of partial differential 
equations for maps $F: U\subset \a \to \b$, $U$ open set of $\a$,
\begin{equation}\label{terngsystem}
  [A_i,\partial_{x^j}F]-[A_j,\partial_{x^i}F] = [[A_i,F],[A_j,F]],\quad 1\le i<j\le n,
   \end{equation}
where $(A_1,\dots,A_n)$ is a basis of regular elements of $\a$ and
$(x^1,\dots,x^n)$ denote the corresponding coordinates. 
This system in known in the literature as the \textit{$G/G_0$-system} 
\cite{BDPT02,TW03} and,
in a different formulation, as the the \textit{curved flat system}
\cite{FP1,FP2}. We show that \eqref{terngsystem} fits into 
our general scheme by proving the following.

\begin{prop}
The tableau $\A$ is involutive, $\Phi$ takes values in 
$B^{0,2}(\A)$ and satisfies \eqref{torsionezero}.
\end{prop}

\begin{proof}
A direct computation shows that the characters of $\A$ are  
$$
 s'_1 = m-n=\dim \b, \quad s'_2=\cdots =s'_n=0.
  $$
Next, observe that an element $F \in \hom(\a,\b)$ 
belongs to $\A^{(1)}$ if and only if
$[F(A'),A'']-[A',F(A'')]=0$, for all $A',A''\in \a$.
On the other hand, from Jacobi identity we get
$[\text{ad}_X(A'),A'']-[A',\text{ad}_X(A'')]=0$, for all $X\in
\p$, $A',A''\in \a$. Therefore, the map $\p \ni X \to \text{ad}_X|_{\a}\in
{\A}^{(1)}$ establishes an isomorphism between ${\A}^{(1)}$ and $\mathbf{P}$. 
This proves the involutiveness of $\A$ by counting dimensions.

For the second statement, note that
$$
 2[[A_1,B],[A_2,B]]=[[B,[A_1,B]],A_2]-[[B,[A_2,B]],A_1],
   $$ 
for $A_1,A_2\in \a$, $B\in \b$, and that $[B,[A_1,B]]\in \m$.
Next, define $S_{|_B}\in C^{1,1}(\A)$ by posing
$S_{|_B}(A): =[B,[A,B]]_{|_{\b}}$. We then have
$2\Phi|_B(A_1,A_2) =\delta^{(1,1)}(S_{|_B})(A_1,A_2)$, for all
$A_1,A_2\in \a$.  This proves that $\Phi|_B\in
B^{0,2}(\A)$, for each $B\in \b$.
From the previous description of ${\A}^{(1)}$,
we get
$$
\Phi_*|_{(A,B)}(A_1+\mathrm{ad}_X(A_1))(A_2,A_3)
   =[[A_2,[X,A_1]],[A_3,B]]+[[A_2,B],[A_3,[X,A_1]].
  $$
Again using the Jacobi identity, we conclude that 
$\Phi$ satisfies \eqref{torsionezero}.

\end{proof}

\subsection{$1+1$ wave maps into Lie groups}

A $1+1$ \textit{wave map} is a smooth map $g:\R^{1,1} \to G$ into a Lie group $G$ 
that satisfies the differential equation
 $$
    (g^{-1}\partial_xg)_y+(g^{-1}\partial_yg)_x=0
   $$
(see \cite{TU}).
If $\theta = g^{-1}dg = Ad x+ Bdy$, then $A, B :\R^{1,1}\to \g$ 
satisfy the first order system
\begin{equation}\label{wavemaps}
  \partial_x B= -[A,B],\quad \partial_y A = [A,B].
   \end{equation} 
Conversely, if $A, B :\R^{1,1}\to \g$ is a solution to \eqref{wavemaps}, then
$\theta = Adx+Bdy$ satisfies the Maurer-Cartan equation and there
exists a unique wave map $g:\R^{1,1}\to G$ such that $g^{-1}dg=\theta$ and
$g(0,0) = e_G$ (the group identity).

  Let $\a =\R^{1,1}$, $\b=\g \oplus \g$ and consider the
tableau 
$$
 \A =\{(X_2dy,X_1dx): X_1,X_2\in \g\}\subset \hom(\a,\b),
   $$ 
which we identify with $\b$.
Next, let 
$$
 \Phi : (A,B) \in \A
   \to - ([A,B],[A,B])^T dx\wedge dy\in \Lambda^2(\a^\ast)\otimes \b.
     $$
It is easily seen that \eqref{wavemaps} amounts to the system
associated with the tableau $\A$ and the non-homogeneous term $\Phi$.
Since there are two independent variables, $\Phi$
automatically satisfies \eqref{torsionezero}.
We need to prove that the tableau is $2$-acyclic and
that $\Phi$ takes values in $B^{0,2}(\A)$.
The first and second prolongation of $\A \simeq \b$ are given, respectively, by
$$
 {\A}^{(1)}=\{Q_{(1)}(\dot{X}_1,\dot{X}_2) =
   (\dot{X}_2dy,\dot{X}_1dx) : \dot{X}_1,\dot{X}_2\in \g\}\subset
    \Lambda^1(\a^\ast)\otimes \A
      $$ 
and
$$
  {\A}^{(2)}=\{Q_{(2)}\, |\, Q_{(2)}(\ddot{X}_1,\ddot{X}_2)=(\ddot{X}_2dy\odot
     dy , \ddot{X}_1dx\odot dx) :\ddot{X}_1,\ddot{X}_2\in \g\}\subset
       S^2(\a^\ast)\otimes \A.
       $$
Therefore,
$\dim {\A}^{(h)} = 2n$, for $h=0,1,2$,
$s_1^0 = n$, $s_2^0 =n$, $s_1^{(1)} = 2n$ and $s_2^{(1)} = 0$. We
then conclude that the involutiveness index of $\A$ is
$k=1$ and that the involutive characters are $s^{(1)}_1=2n$, $s^{(1)}_2=0$.
Since ${\A}^{(1)}$ is involutive, to prove that $\A$ is $2$-acyclic
we only need to prove that $H^{1,2}(\A) = (0)$.
In the case of two independent variables, this amounts
to prove the surjectivity of $\delta^{2,1}$. 
For this, it is sufficient to observe that 
$$
 C^{2,1}(\A)=\{\mu=(\mathcal{X}_2 \otimes dy, \mathcal{X}_1\otimes dx) :
\mathcal{X}_i=X_{i1}dx+X_{i2}dy\in \Lambda^1(\a^\ast)\otimes \g\}
   $$ 
and that $\delta^{(2,1)}(\mu) = (X_{21},X_{12})dx\wedge dy$.
Finally, we show that $\Phi$ takes values in $B^{0,2}(\A)$ by proving 
the vanishing of $H^{0,2}(\A)$, which in turn amounts to proving that 
$\delta^{1,1}$ is onto.
The claim follows by observing that
$$
  C^{1,1}(\A)=\{\sigma = (\mathcal{X}_2,\mathcal{X}_1):
   \mathcal{X}_i=X_{i1}dx+X_{i2}dy\in \Lambda^1(\a^\ast)\otimes \g\}
     $$ 
and that $\delta^{1,1}(\sigma) = (X_{2,1},-X_{12})dx\wedge dy$.

\section{The Cauchy problem}\label{s:cauchy}

\subsection{Statement of the result}

Retaining the notation of Sections \ref{s:basic} and \ref{s:ideals},
let $\mathcal{A}=(A_1,\dots,A_n)$ be a basis of $\a$ and 
let $\mathcal{Q}^{h} = (Q^{(h)}_1,\dots,Q^{(h)}_{r_h})$ be a basis
of the prolongation ${\A}^{(h)}$ $(h= 0,1,\dots,k)$, where $k$ is the 
least integer such that the prolongation ${\A^{(k)}}$ is involutive.
Next, let
$\mathcal{Q}^k=(Q^{(k)}_{[1],1},\dots, Q^{(k)}_{[1],s_1},
\dots,Q^{(k)}_{[\nu],1},\dots,Q^{(k)}_{[\nu],s_{\nu}})$ be
a normal basis of ${\A}^{(k)}$ and suppose that its elements be in Guillemin's
normal form with respect to 
$\mathcal{A}$ and $\mathcal{Q}^{k-1}$ (cf. Section \ref{ss:guillemin}). 

We call $(\mathcal{A},\mathcal{Q}^{0},\dots,\mathcal{Q}^{k})$ a
\textit{regular basis} of $M_{(k)}$. Let
${\A}^{(k)}={\A}^{(k)}_{[1]}\oplus \cdots\oplus
{\A}^{(k)}_{[\nu]}$, where
$\mathbf{A}^{(k)}_{[p]}=\span(Q^{(k)}_{[p],1},\dots,Q^{(k)}_{[p],s_p})$,
$p=1,\dots,\nu$. Consider the filtration
$$
 \a_1=\span(A_1)\subset
  \a_2=\span(A_1,A_2)\subset\cdots\subset
   \a_{\nu}=\span(A_1,\dots,A_{\nu})\subset \a ;
     $$ 
and denote by $\pi_{\rho} : \a\to \a_{\rho}$ the projection of $\a$ on 
$\a_{\rho}$ with respect to the decomposition
$\a=\a_{\rho}\oplus\span(A_{\rho+1},\dots,A_n)$.

\begin{defn}
A \textit{set of Cauchy data} consists of an ordered set
$$
 (U,x_0, P_{(0)},\dots,P_{(k-1)},P^{[1]}_{(k)},\dots,P^{[\nu]}_{(k)}),
   $$
where $U\subset \a$ is an open set of  
$x_0\in \a$, $P_{(0)}\in {\A}^{(0)},\dots,P_{(k-1)}\in
{\A}^{(k-1)}$ are constants, and  
$P^{[\rho]}_{(k)}: \pi_{\rho}(U)\subset \a_{\rho}\to
{\A}^{(k)}_{[\rho]}$ ($\rho =1,\dots,\nu$) are analytic maps.
\end{defn}

\begin{thm}\label{cauchypbm}
Let $\A \subset \hom(\a,\b)$ be a $2$-acyclic tableau,
$k$ the least integer such that the prolongation
$\A^{(k)}$ is involutive and 
$s_1^{(k)},\dots,s^{(k)}_{\nu}$ the corresponding principal characters. 
Let $\Phi:\a\oplus \A\to \b \otimes \Lambda^2(\a^\ast) \subset B^{0,2}$ 
be an analytic function satisfying \eqref{torsionezero} and
$(M_{(k)},\IG_{(k)},\Omega)$ the involutive prolongation
of the differential system $(M,\IG, \Omega)$ defined by $\A$ and $\Phi$ (cf. Theorem
\ref{theoremA}). Let 
$$
  (U,x_0,P_{(0)},\dots,P_{(k-1)},P^{[1]}_{(k)},\dots,P^{[\nu]}_{(k)})
   $$ 
be a set of analytic Cauchy data. Then there exist an open set
$\widetilde{U}\subset U$, $x_0\in\widetilde{U}$, and an analytic map
$Q=(Q_{(0)},\dots,Q_{(k)}):U\to {\A}^{[k]}$ such that:

\begin{enumerate}
\item $\widetilde{U} \ni x \mapsto (x,Q(x))\in M_{(k)}$ is a $K$-regular
integral manifold of $(\IG_{(k)},\Omega)$;

\item $Q$ satisfies the initial conditions
$Q_{(h)}(x_0)=P_{(h)}$, $h=0,\dots,k-1$ and $Q_{(k)}^{[\rho]}\circ
\pi_{\rho} = P_{(k)}^{[\rho]}$, for all $\rho=1,\dots,\nu$;

\item $Q$ is unique, in the sense that any other analytic map
with the same properties coincides with $Q$ in a neighborhood
of $x_0$.

\end{enumerate}
\end{thm}

\begin{remark}
If the \textit{characteristic variety} of 
${\A}^{(k)}$ is hyperbolic (cf. \cite{YA}), the result holds also
in the case of ${C}^{\infty}$ data.
\end{remark}

\subsection{Proof of Theorem \ref{cauchypbm}}

The proof will be divided in several steps.
Consider the vector bundles defined by
$$
 J^{\perp}=\ker(\beta_{(0)},\dots,\beta_{(k-1)},dx^1,\dots,dx^n)
  \subset I^{\perp}=\ker(\beta_{(0)},\dots,\beta_{(k-1)})\subset
    TM_{(k)}.
     $$ 
Let $G_h(I,J) \to M_{(k)}$ be the Grassmann bundle of
$h$-dimensional subspaces of $I^{\perp}$ transverse to
$J^{\perp}$ and denote by $V_h(I,J)\subset G_h(I,J)$ the set of 
$h$-dimensional integral elements of $(\IG_{(k)},\Omega)$.
Let
$$
 (P,\xi)\in M_{(k)}\times \a\to V(\xi)_P=\xi
  +\sum_{h=1}^{k} (S_{(h)}|_P+Q_h)(\xi)
  $$
and consider the trivialization of $I^{\perp}$ given by
$$
 (P,\xi,\mathcal{X})\in M_{(k)}\times (\a\oplus
   \mathbf{A}^{(k)})\to V(\xi)_P+\mathcal{X}\in I^{\perp},
    $$ 
for $P=(x,Q_{(0)},\dots,Q_{(k)})\in M_{(k)}$, $\xi\in \a$,
$\mathcal{X}\in {\A}^{(k)}$, and the identifications
$$
 G_h(I,J)\simeq M\times G_h(\a\oplus
   {\A}^{(k)},{\A}^{(k)}), 
      \quad G_h(I^{\perp}/J^{\perp})\simeq M\times G_h(\a),
         $$ 
where $G_h(\a\oplus {\A}^{(k)},{\A}^{(k)})$ denotes the Grassmannian of
$h$-dimensional subspaces in $\a\oplus {\A}^{(k)}$ transverse to ${\A}^{(k)}$.
Further, $G_h(\a\oplus {\A}^{(k)},{\A}^{(k)})$ may be identified
with the fiber bundle $\hom(\mathcal{T}_h,{\A}^{(k)})$ via the 
isomorphism
$$
  (\a_h,F)\in \hom(\mathcal{T}_h,{\A}^{(k)})\to 
   \mathcal{A}_h(\a_h,F):=\{\xi+F(\xi) : \xi \in \a_h\}\in G_h(\a\oplus
       {\A}^{(k)},{\A}^{(k)}),
         $$ 
where $\mathcal{T}_h\to G_h(\a)$ denotes the tautological bundle.
Therefore, we can think of $V_h(I,J)$ as an analytic subset
of $M_{(k)}\times \hom(\mathcal{T}_h,{\A}^{(k)})$. 
Next, let $G^\sharp_h(\a)\subset G_h(\a)$ be the Zariski open of
$h$-dimensional generic subspaces with respect to the tableau
${\A}^{(k)}$ and let $V_h^{\sharp}(I,J)\subset V_h(I,J)$
the Zariski open of integral elements that project on $h$-dimensional
generic subspaces.

\begin{lemma}
The projection $V_h^{\sharp}(I,J)\to G^\sharp_h(\a)$  
is an (analytic) affine bundle, whose fibers have dimension
 $s_1+2s_2+\cdots+(h-1)s_h+h(s_h+\cdots+s_n)$. 
\end{lemma}

\begin{proof}
The relative tableau ${\A}^{(p)}\| _{\a _h}:=\{Q|_{\a_h} \, |\,
Q\in {\A}^{(p)}\}$ of a generic subspace $\a_h\subset
\a$ has dimension $s_1+\cdots+s_p$. Therefore,
$$
  {\mathcal{A}}_h^{(p)}=\{(\a_h,{Q_{(p)}}_{|_{\a_h}}) \,:\, 
     \a_h \in G^\sharp_h(\a), Q_{(p)}\in
      {\A}^{(p)}\}\to G^\sharp_h(\a)
       $$ 
is a vector bundle of rank $s_1+\cdots+s_p$. The involutiveness of
${\A}^{(k)}$ implies that the morphism $\hat{\delta}:
\hom(\mathcal{T}^{\sharp}_h,{\mathcal{A}}^{(k)}_h)\to
\Lambda^2(\mathcal{T}^{\sharp}_h)\otimes {\mathcal{A}}^{(k-1)}_h$
defined by
$$
 \hat{\delta}(F)(\xi_1,\xi_2)=F(\xi_1)(\xi_2)-F(\xi_2)(\xi_1),
   $$
for $F\in \hom(\a_h,\mathcal{A}_h^{(k)}|_{\a_h})$, and
$\xi_1,\xi_2\in \a_h$, has constant rank.

 The kernel of $\hat{\delta}$, which we denote by 
$({\mathcal{A}}_h^{(k+1)})^{(1)}$, is a subbundle, whose fiber above
$\a_h$ is the first prolongation of ${\A}^{(k)}_{\|_{\a_h}}$. 
To compute the dimension of the fibers it suffices to observe that
${\A}^{(k)}_{\| _{\a_h}}$ is involutive, from which it follows that 
$\dim({\A}^{(k)}_{\|_{\a_h}})^{(1)}=s_1+2s_2+\cdots+hs_h$. 

Let $\rho :\hom(\mathcal{T}_h^{\sharp},{\A}^{(k)})\to
\hom(\mathcal{T}^{\sharp}_h,{\mathcal{A}}^{(k)}_h)$ be the restriction morphism 
defined by 
$$
 \rho(F)(\xi)=F(\xi)|_{\a_h}, \text{ for } F\in
   \hom(\a_h, {\A}^{(k)}), \, \xi\in \a_h, \,
    \a_h\in G^\sharp_h(\a).
      $$ 
Let
$\Psi:\hom(\mathcal{T}^{\sharp}_h,\mathbf{A}^{(k)})\to
\hom(\mathcal{T}^{\sharp}_h,\mathcal{A}_h^{(k)})/(\mathcal{A}_h^{(k)})^{(1)}$
be the composition of $\rho$ with the projection on the quotient and let
$\hat{S}_{(k+1)}:M_{(k)}\times G_h(\a)\to M_{(k)}\times
\text{Hom}(\mathcal{T}_h,\mathbf{A}^{(k)})$ be the section defined by 
$(P,\a_h)\to (P,S_{(k+1)}|_P\|_{\a_h})$. 
The proof of the lemma follows by observing that
$$
 V^{\sharp}(I,J)=\{(P,F)\in M_{(k)}\times
   \hom(\mathcal{T}^{\sharp}_h,\mathbf{A}^{(k)}) \, : \,
    \Psi(F-\hat{S}_{(k+1)}|_P)=0\}
     $$ 
and that
$\mathrm{rank}(\ker(\Psi))=s_1+2s_2+\cdots+(h-1)s_{h-1}+h(s_h+\cdots+s_n)$.
\end{proof}

\begin{lemma}
The generic integral flags are $K$-regular. 
\end{lemma}

\begin{proof}
For any $P\in M_{(k)}$ and for any $\mathcal{A}_h\in G_h^{\sharp}(\a\oplus
\mathbf{A}^{(k)},\mathbf{A}^{(k)})$ consider the \textit{polar space}
$$
  \mathcal{H}(P,\mathcal{A}_h): =
 \{\xi_1+\mathcal{X}_1\in \a\oplus \mathbf{A}^{(k)} \, : \,
    (\pi_{(k)}\dot{\wedge}\eta)(\xi_1+\mathcal{X}_1, \xi +
\mathcal{X})=0,\text{ for all }\xi+\mathcal{X}\in \mathcal{A}_h \}.
$$
Since the projection of $\mathcal{A}_h$ on $\a$ is a generic subspace,
we have that 
 $$
  \dim(\mathcal{H}(P,\mathcal{A}_h))=n+(s_{h+1}+\dots+s_{\nu}).
    $$
If $\mathcal{A}_0=(0)\subset \mathcal{A}_1\subset\cdots\subset
\mathcal{A}_n\subset \overline{\a}\oplus \mathbf{A}^{(k+1)}$ is an integral flag 
in $P\in M_{(k)}$  and if the $\mathcal{A}_h$ are generic, then 
$$
  \dim\mathcal{H}(P,\mathcal{A}_h) = n+(s_{h+1}+\cdots+s_{\nu}),
  $$
for $h=0,\dots,n$. From this identities and the preceding lemma 
it follows that
$(\mathcal{A}_0,\dots,\mathcal{A}_n)$ is an integral flag which is
$K$-regular.
\end{proof}

\begin{lemma}
let  $(\mathcal{A},\mathcal{Q}^{0},\dots,\mathcal{Q}^{k})$ be
a regular basis and $\mathcal{A}_{h+1}\subset
\a\oplus \mathbf{A}^{(k)}_{[h+1]}\oplus \cdots\oplus
\mathbf{A}^{(k)}_{[\nu]}$ an $(h+1)$-dimensional subspace.
Assume that $\mathcal{A}_{h+1}$ projects onto $\a_{h+1}\subset
\a$ and set
$$
 \Sigma(\mathcal{A}_{h+1}) :=\mathbf{A}^{(0)}\oplus\cdots\oplus
  \mathbf{A}^{(k-1)}\oplus (\mathcal{A}_{h+1}\oplus
   \mathbf{A}^{(k)}_{[1]}\oplus\cdots\oplus\mathbf{A}^{(k)}_{[h]}).
     $$ 
If $\mathcal{A}_h$ is an $h$-dimensional integral element at
$P\in M_{(k)}$ contained in $\Sigma(\mathcal{A}_{h+1})$ and if
$\mathcal{A}_h$ projects onto $\a_h\subset \a_{h+1}$, then
$\dim(\mathcal{H}(P,\mathcal{A}_{h})\cap \Sigma(\mathcal{A}_{h+1}))=h+1$.
\end{lemma}

\begin{proof}
First observe that $\mathcal{H}(P,\mathcal{A}_{h})\cap \Sigma(\mathcal{A}_{h+1})=
\mathcal{H}(P,\mathcal{A}_{h})\cap (\Sigma(\mathcal{A}_{h+1})\cap
{I^{\perp}}_{|_P})$, and that the elements of $\Sigma(\mathcal{A}_{h+1})\cap
{I^{\perp}}_{|_P}$ are of the form
$V_P(\xi)+\mathcal{X}(\xi)+\mathcal{S}$, where $\xi \in \a_{h+1}$,
$\mathcal{X}\in \hom(\a_{h+1},
\mathbf{A}^{(k)}_{[h+1]}\oplus\cdots\oplus
\mathbf{A}^{(k)}_{[\nu]})$ and $\mathcal{S}\in
\mathbf{A}^{(k)}_{[1]}\oplus \cdots\oplus \mathbf{A}^{(k)}_{[h]}$. 
Next, 
the elements of $\mathcal{A}_h$ are of the form
$V_P(\xi)+\mathcal{X}(\xi)+\mathcal{Y}(\xi)$, where $\xi\in \a_{h}$
and 
$\mathcal{Y}\in \hom(\a_h, {\A}^{(k)}_{[1]}\oplus
\cdots\oplus \mathbf{A}^{(k)}_{[h]})$.
If $\r\subset \a_{h+1}$ is a $1$-space such that 
$\a_{h+1}=\r\oplus \a_h$, the elements of
$\Sigma(\mathcal{A}_{h+1})\cap {I^{\perp}}_{|_P}$ are the form
$$
  Y(\xi,\xi',\mathcal{S})=(V_P(\xi)+\mathcal{X}(\xi)+
     \mathcal{Y}(\chi))+(V_P(\xi')+\mathcal{X}(\xi'))+\mathcal{S},
       $$
where $\xi\in \a_h$, $\xi'\in \r$ and $\mathcal{S}\in
\mathbf{A}^{(k+1)}_{[1]}\oplus\cdots \oplus\mathbf{A}^{(k+1)}_{[h]}$. 

If we impose that
$Y(\xi,\xi',\mathcal{S})$ satisfy the polar equations, then 
(using the fact that the elements of the basis $\mathcal{Q}^{k}$
 are in normal form) we find that $\mathcal{S}$ is determined by
 $\xi$ and then
$\dim(\mathcal{H}(P,\mathcal{A}_{h})\cap \Sigma(\mathcal{A}_{h+1}))=h+1$.
\end{proof}

Let $(x,Q_{(0)},\dots,Q_{(k)},Q_{(k)}^{[1]},\dots,Q_{(k)}^{[\nu]})$ be the
coordinates with respect to the regular basis. Consider a set $(U,x_0,
P_{(0)},\dots,P_{(k-1)},P^{[1]}_{(k)},\dots,P^{[\nu]}_{(k)})$ of
analytic Cauchy data and the filtration of submanifolds
$$
 \mathcal{N}_0\subset
  \mathcal{N}_1\subset\cdots\subset \mathcal{N}_{\nu}\subset
   \mathcal{N}_{\nu+1}\subset \cdots\subset \mathcal{N}_n=M
     $$ 
defined as follows:
\begin{itemize}
\item $\mathcal{N}_0$ is the point of coordinates $x=x_0$,
$Q^{(h)}=P_{(h)}$,
$Q^{[p]}_{(k)}=P^{[p]}_{(k)}(x_0^1,\dots,x_0^{p})$, for $h=1,\dots,k-1$ and
$p=1,\dots,\nu$; 

\item  for $r=1,\dots,\nu$, $\mathcal{N}_r$ is defined
by the equations $x^{r+1} = x_0^{r+1},\dots,x^{n}=x_0^{n}$,
$Q^{[r+p]}_{(k)}=P^{[r+p]}_{(k)}(x^1,\dots,x^r,x_0^{r+1},\dots,x_0^{\nu})$,
for $p=1,\dots,\nu-r$; 

\item for $r=\nu+1,\dots,n$, $\mathcal{N}_r$ 
is defined by $x^{r+1}=x_0^{r+1},\dots,x^{n}=x_0^{n}$. 
\end{itemize}

\noindent From the preceding lemma it follows that if $\mathcal{A}_h\subset
T_P(\mathcal{N}_{h+1})$ is an $h$-dimensional integral element,
then $\dim(\mathcal{H}(P,\mathcal{A}_h)\cap
T_P(\mathcal{N}_{h+1}))= h+1$. 
By applying the Cartan-K\"ahler theorem we prove the existence of
a unique filtration of $K$-regular
integral submanifolds $P_1\subset \cdots\subset
P_{\nu}\subset P_{\nu +1}\subset \cdots\subset P_n$,
$\dim P_h =h$, such that $\mathcal{N}_{h-1}\subset P_h\subset
\mathcal{N}_h$, for $h=1,\dots,n$. 

Let $\pi_1 : P_n\to \a$ and 
$\pi_2  : P_n\to {\A}^{[k]}$ be the restrictions of the respective
projections of $M_{(k)}$. 
Since $\pi_1$ has maximal rank we may assume that the image
be an open neighborhood $\widetilde{U}$ of $x_0$ and that
$\pi_1$ be invertible. Thus $Q = \pi_2\circ
(\pi_1)^{-1}:\widetilde{U}\to \mathbf{A}^{[k]}$ satisfies
the requested properties.
The uniqueness of $Q$ is a consequence of the uniqueness
of the $K$-regular integral submanifolds
$P_1\subset \cdots\subset P_n$ such that $\mathcal{N}_{h-1}\subset
P_h\subset \mathcal{N}_h$, $h=1,\dots,n$. This proves the theorem.

\bibliographystyle{amsalpha}

\end{document}